\documentclass[12pt]{article}
\usepackage{amsmath,amsfonts,latexsym,amssymb}

\numberwithin{equation}{section}

\newcommand{\C}{{\mathbb C}}

\newcommand{\Z}{{\mathbb Z}}

\begin{document}
\begin{center}
{\Large{\bf On Representations of Toroidal Lie Algebras}} \\ [2cm]
{\bf S. Eswara Rao} \\
{\bf School of Mathematics} \\
{\bf Tata Institute of Fundamental Research} \\
{\bf Homi Bhabha Road} \\
{\bf Mumbai - 400 005}
 {\bf India} \\ [5mm]
{\bf email: senapati@math.tifr.res.in}
\end{center}

\section*{$0 \ \ \S$ Introduction:}

\hspace*{5mm} I gave two talks on some of 
my results on toroidal Lie algebras in the conference {\bf Functional Analysis VIII} held in Dubrovnik, 
Croatia in June 2003.
  But these results have been submitted to research
journals and will appear soon.  So I decided to write an expository 
article on Toroidal Lie Algebras covering my results for the proceedings
of Functional Analysis VIII.

I will introduce the definitions of toroidal Lie algebras and also the 
generalized Virasoro algebra in the introduction below.  In the main
body of the article I will state most of the recent results on 
representations of toroidal Lie algebras with finite dimensional weight
spaces.

Toroidal Lie algebras are $n$-variable generalizations of the well known 
affine Kac-Moody Lie algebras.  The affine Kac-Moody Lie algebra,
which is the universal central extension of Loop algebra, has a very
rich theory of highest weight modules and some of their characters admit
 modular properties.  The level one highest weight integrable modules has been
constructed explicitly on the Fock space through the use of vertex operators
(see [FK]; for an introduction to vertex operator algebras, see [LL]).
In such cases the characters can be computed easily and using Weyl - Kac
character formula interesting identities have been obtained.

There is another well known Lie algebra called Virasoro Lie algebra
which is the universal central extension of the Lie algebra of 
diffeomorphisms of one dimensional torus.  The Virasoro Lie algebra acts on any (except when the level is negative of dual coxeter number) highest weight
module of the affine Lie algebra through the use of famous Sugawara operators.  So the semi-direct product of affine Lie algebra and Virasoro Lie
algebra with common center has turned out to be a very important Lie algebra
in mathematics and also in physics.  This Lie algebra plays a very significant role in Conformal Field Theory.  The connection to physics is explained
in great detail by Di Francesco, Mathieu and Senechal in the book on Conformal Field Theory [DMS].

Toroidal Lie algebras are born out of an attempt to generalize the above classical theory.  We will first define toroidal Lie algebra and explain the main
results obtained.

Let ${\cal G}$ be a simple finite dimensional Lie algebra and let $A=\C[t_1^{\pm 1}, \cdots, t_n^{\pm 1}]$ be a Laurent polynomial ring in $n$ 
commutating variables.  Then ${\cal G} \otimes A$ has a natural structure 
of Lie algebra and the universal central extension is what we call
{\it Toroidal Lie algebra}.  For the first time a large class of 
representations are constructed in [EMY]  and [EM] through the use 
of vertex operators on the Fock space.  Thus giving a realization
of toroidal Lie algebras.

Affine Lie algebras are precisely the case when $n=1$.  One significant
difference for toroidal Lie algebra is that the universal center for
$n \geq 2$ is infinite dimensional where as it is one dimensional for 
affine Lie algebras.  Toroidal Lie algebras are naturally $\Z^n$-graded and
 there is no natural decomposition into positive and negative spaces.
Thus the standard highest weight module theory does not apply for 
toroidal case.  Also the infinite center does not act as scalars on
irreducible module but as invertible operators.  So the study of
modules does not follow the standard methods of affine modules and needs
 completely different techniques.

Surprisingly, a lot of information can be obtained on an irreducible module 
for toroidal Lie algebra when we assume all {\it weight spaces are finite
dimensional.}  For example the part of the center which acts non-trivially
on an irreducible module will have a definite shape.  In the sense that
if we mod out some arbitrary part of the center and look for an
irreducible module with finite dimensional weight spaces, they may not exist.
For example if we mod out all the non-zero degree center then module
does not exist with finite dimensional weight spaces and with non-zero
central (zero degree) action.  Further it is proved that the ratio of 
standard zero degree central operators on an irreducible module with
finite dimensional weight spaces is rational. See ref. [E8].

The most striking result is the classification of irreducible integrable
modules for toroidal Lie algebras with finite dimensional weight spaces.
But the classification is disappointing as the modules are tensor product 
of evaluation modules of affine Lie algebras.  In particular the 
characters are product of affine characters.  See references [E4,E5 and
Y].

It is worth noting that integrable modules are not completely reducible and
such modules are constructed by Chari and Le Thang [CL].  So the 
classification of integrable modules not necessarily irreducible may be 
interesting and their characters being Weyl group invariant may hold 
some significance.

The next natural issue is to generalize the Virasoro Lie algebra. Thus we consider the Lie algebra of diffeomorphisms on $n$ dimensional torus which is known to be isomorphic to $DerA$ the derivation algebra of $A$. 
Several attempts have been made by physicists to give a  Fock space 
 representation to Der$A$ or to its extension (see [FR]).  They all failed to produce
any interesting results due to lack of proper definitions of ``normal
ordering'' among other things.  At this juncture an interesting result
has come out in [RSS] which says that Der$A$ has no non-trivial central extension.

Let us go back to the Vertex construction of toroidal Lie algebra in [EM]
where operators are defined for Der$A$ on the Fock space generalizing
the Sugawara construction.  But the corresponding extension for Der$A$
turned out to be very wild (certainly non-central) and not tractable [E7].
In the process an interesting abelian extension for Der$A$ has been
created in [EM] and the abelian part is exactly  the center of the
toroidal Lie algebra.  So the semi-direct sum of toroidal Lie algebra and the Der$A$ with common extension has emerged as an interesting object of study 
which we will now define.

First we define a toroidal Lie algebra see [MY] and [Ka].  Let ${\cal G}$
be a simple finite dimensional Lie algebra over complex numbers.  Let $<,>$ 
be a ${\cal G}$-invariant symmetric non-degenerate bilinear form on ${\cal G}$.
Fix a positive integer $n$ and recall $A=\C[t_1^{\pm 1}, \cdots, t_n^{\pm 1}]$.
Let $\Omega_A$ be the module of differentials which can be defined as
vector space spanned by the symbols  $\{t^rK_i, r \in \Z^n, \ 1 
\leq i \leq n \}$. Let $d_A$ be the subspace spanned by $\sum r_i t^rK_i$ 
and consider $\Omega_A/d_A$.  Let $K(u,r)= \sum u_i t^rK_i, u=
(u_1, \cdots, u_n) \in \C^n$.  Then the toroidal Lie algebra is
$\tau = {\cal G} \otimes A \oplus \Omega_A/d_A$
with Lie bracket

\begin{enumerate}
\item[(1)] $[X \otimes t^r, Y \otimes t^s] = [X,Y] \otimes t^{r+s}+<X,Y> K (r,r+s)$.
\item[(2)] $K (u,r)$ is central. \\
Here $X \in {\cal G}$ and $t^r=t_1^{r_1} \cdots t_n^{r_n} \in A$.  
\end{enumerate}

Der$A$ has $\{t^r t_i \frac{d}{dt_i}, r \in \Z^n, 1 \leq i \leq n\}$ for a basis.
Write $D(u,r) = \sum u_i t^r t_i \frac{d}{dt_i}, u \in  \C^n$.  The
abelian extension $\widehat{DerA}$ of Der$A$ is as follows [EM].
$$\widehat{DerA} = \Omega_A/d_A \oplus DerA.$$
The Lie bracket is given by,
\begin{enumerate}

\item[(3)] $[D(u,r), D(v,s)]=D(w,r+s) - (u,s) (v,r) K(r,r+s)$, \\
where $w=(u,s) v-(v,r)u$ 
\item[(4)] $[D(u,r), K(v,s)] = (u,s) K (v,r+s)+ (u,v)K (r,r+s)$
\item[(5)] $[K (u,r), K(v,s)]=0$ 
\end{enumerate}
where $(,)$ is the standard bilinear form on $\C^n$.  We will now define the 
semi-direct product of toroidal Lie algebra and Der$A$ with common extension
which we call full toroidal Lie algebra.
$$\widehat{\tau} = {\cal G} \otimes A  \oplus \Omega_A/d_A \oplus DerA$$
\begin{enumerate}
\item[(6)] $[D (u,r), X \otimes t^s ] =(u,s) X \otimes t^{r+s}$ \\
(1) to (6) defines the Lie structure on $\widehat{\tau}$.
\end{enumerate}
The first natural question is, does there exists a realization for $\widehat{\tau}$, in other words can we construct a natural module for $\widehat{\tau}$.
Several attempts have been made in [EM], [BB] and [BBS].  Eventually Yuly 
Billig in [B3] has succeeded in constructing a module through the use of
Vertex operator algebras.  The irreducible integrable modules for
$\widehat{\tau}$ with finite dimensional weight spaces has been classified
recently see references [JM2] and [EJ].

\paragraph*{$\S$ 1. Representation of Toroidal Lie algebras}

\paragraph*{(1.1)}  Let ${\cal G}$ be simple finite dimensional Lie algebra over complex numbers.  For any positive integer $n$, let $A=A_n =\C [t_1^{\pm 1},
\cdots, t_n^{\pm 1}]$ be Laurent polynomial ring in  $n$ commuting
variables.  The universal central extension of ${\cal G} \otimes A$ is called
toroidal Lie algebra $\tau$.  $\tau$ can be constructed more explicitly.

For $r= (r_1, \cdots, r_n) \in \Z^n$ let $t^r =t_1^{r_1} t_2^{r_2} \cdots
t_n^{r_n} \in A$.  Let $\Omega_A$ be the vector space spanned by the symbols
$\{t^r K_i, r \in \Z^n, 1 \leq i \leq n\}$.  Let $d_A$ be the subspace spanned by $\{\sum r_i t^r K_i, r \in \Z^n \}$.

Let $\tau = {\cal G} \otimes A \oplus \Omega_A/d_A$ and the Lie structure
 is given by
\begin{enumerate}
\item[(1)] $[X \otimes t^r, Y \otimes t^s ]= [X,Y] \otimes t^{r+s} +<X,Y> d 
(t^r)t^s$ where $X,Y \in {\cal G}, r,s \in \Z^n$ and $<,>$ is a ${\cal G}$-invariant symmetric bilinear form on ${\cal G}$.  Further $d(t^r)t^s = \sum r_i t^{r+s}K_i$.
\item[(2)] $\Omega_A /d_A$ is central in $\tau$.
\end{enumerate}

\paragraph*{Theorem} ([MY], [Ka]) $\tau$ is the universal central extension
of ${\cal G} \otimes A$.

Clearly $\tau$ is naturally graded by $\Z^n$.  To reflect this fact let $D$
be the linear span of degree derivations $d_1, \cdots, d_n$ definite as
$$[d_i, X \otimes t^r]= r_i X  \otimes t^{r}$$
$$[d_i, t^r K_j]= r_i t^r K_j$$
$$[d_i, d_j] =0$$
Let $\tilde{\tau} = \tau \oplus D$ and is also called a toroidal Lie algebra.
It is easy to see that
$$\dim (\Omega_A/d_A)_r =n-1 \ {\rm if} \ r \neq 0.$$
$$\dim (\Omega_A / d_A)_0 =n.$$
Thus dimension of $\Omega_A/d_A$ is infinite if $n \geq 2$ and equals to one if $n=1$.  In fact from the definition it follows that $\tilde{\tau}$ is
affine Kac-Moody Lie algebra for $n=1$.

Thus toroidal Lie algebras are $n$ variable generalization of affine 
Kac-Moody Lie algebras.

\paragraph*{(1.2)} Root space decomposition for $\tilde{\tau}$ [E5].
Fix a Cartan subalgebra $h$ of ${\cal G}$ .  Let $\triangle$ be the root system of ${\cal G}$ and let $\alpha_1, \cdots \alpha_d$ be set of simple roots.
Let $\alpha_1^{\vee}, \cdots, \alpha_d^{\vee}$ be simple co-roots.  Let ${\cal G} = \bigoplus_{\alpha \in \triangle} {\cal G}_{\alpha} \oplus h$ be the 
root space decomposition of ${\cal G}$.  We will now describe the root space
 decomposition of $\tilde{\tau}$ with respect to the Cartan $\tilde{h}=
h \oplus \sum \C K_i \oplus D$.  Let $\delta_1, \delta_2, \cdots \delta_n$
be in $\tilde{h}^*$ such that $\delta_i (d_j) = \delta_{ij}, \delta_i 
(\alpha^{\vee}_j) =0$ and $\delta_i (K_j) =0$.  Let $\Lambda_1, \cdots 
\Lambda_n$ be in $\tilde{h}^*$ such that $\Lambda_i (\alpha_j^{\vee})=0, 
\Lambda_i (d_j)=0$ and $\Lambda_i (K_j)= \delta_{ij}$.  Then it is clear 
that $\alpha_1, \cdots \alpha_d, \delta_1, \cdots \delta_n, \Lambda_1, 
\cdots , \Lambda_n$ is a basis of $\tilde{h}^*$ for dimension reasons.  It is now easy to give
a non-degenerate (non-positive definite) symmetric bilinear form on $\tilde{h}^*$ such that the corresponding matrix with respect to the basis $\alpha_1,
\cdots, \alpha_d, \delta_1, \cdots \delta_n, \Lambda_1, \cdots , \Lambda_n$
is
$$
\begin{bmatrix}
 AD,& 0,& 0 \\
0,& 0,& I \\
0,& I,&  0
\end{bmatrix}
$$

Here $A$ is the finite Cartan matrix and $D$ is diagonal so that $AD$ is 
symmetric.

For $r \in \Z^n$ let $\delta_r = \sum r_i \delta_i$ and note that  $(\delta_r,
\delta_s)=0$.  The $\delta_r$'s are called null roots.
Let $\tilde{\triangle} =\{ \alpha + \delta_r, \delta_s \mid \alpha \in \triangle, r, s \in \Z^n \}$.  Then $\tilde{\tau}$ has the following root space
decomposition

$\tilde{\tau} = \displaystyle{\bigoplus_{\gamma \in \tilde{\triangle}}} \ 
\tilde{\tau}_\gamma$ where \\
if $\gamma = \alpha + \delta_r, \alpha \in \triangle$ then $\tilde{\tau}_\gamma= {\cal G}_{\alpha} \otimes t^r$ \
if $\gamma =\delta_r \neq 0$ then $\tilde{\tau}_r = h \otimes t^r$ \\
if $\gamma =0$ then $\tilde{\tau}_0 = \tilde{h}$.

Let $\alpha_{d+j} = \beta - \delta_j, 1 \leq j \leq n$ which  should be 
thought of as additional simple roots.  It is  easy to see that every root is 
linear combinations of $\alpha_1, \cdots, \alpha_{d+n}$ but not non-negative
or non-positive linear combinations of $\alpha$'s.  Thus there is no natural 
positive or negative root spaces.  So the standard highest weight module 
theory does not go through for toroidal Lie algebra.

If we compute the matrix $(\alpha_i, \alpha_j)_{1 \leq i, j \leq d+n}$ one see
that it is a generalized intersection matrix (GIM).  Thus toroidal Lie 
algebras are quotients of GIM algebras as defined by Slodowy.  In fact
they 
are GIM algebras in the symmetric case and it is open problem in the 
non-symmetric case.

We will now describe some of the known results on toroidal Lie algebras.  
The first is the shape of the possible non-zero center on an arbitrary 
irreducible module.  The second one is the classification of irreducible 
integrable modules.  The third one is the sketch of the construction of 
faithful representation for toroidal Lie algebra through the use of 
vertex operators.

\paragraph*{(1.3)}  Shape of the non-zero center of an irreducible 
representation [E8].  We fix an irreducible representation $V$ for the 
toroidal Lie algebra with finite dimensional weight spaces.  An operator 
$z:V \to V$ is called central operator of degree $m \in \Z^n$ if $z$ 
commutes with $\tau$ action and $d_i z - zd_i=m_i z$.  For example $t^m K_i$ 
is a central operator of degree $m$ and $t^n K_i t^m K_j$ is a central operator of degree $m+n$.
\paragraph*{Lemma (a)}  Suppose a central operator is non-zero at one vector
then it is non-zero at every vector. \\
(b) A non-zero central operator of degree $m$ admits an inverse which is a 
central operator of degree $-m$.  \\
(c)  Suppose $z_1$ and $z_2$ are non-zero central operators of degree $m$ 
then $z_1 = \lambda z_2$ for some non-zero scalar $\lambda$.

Let $L=\{r \in \Z^n \mid t^r K_i $ is non-zero on $V$ for some $i\}$.
Let $S$ be the group generated  by $L$.  We first note that for $m \in S$ there is a non-zero central operator of degree $m$.

Let rank of $S=k$. \\
Before we explain our results we will first explain the notion of change
of co-ordinates.
\paragraph*{Choice of co-ordinates: } Let $B \in GL(n, \Z)$ and let $e_1, 
\cdots e_n$ be the standard basis of $\Z^n$.

Let $s_i =t^{B (e_i)}$. Then the subalgebra generated by $s_i^{\pm 1}$ 
which is $\C [s_1^{\pm 1}, \cdots s_n^{\pm 1}]$ and is equal to 
$\C [t_1^{\pm 1}, \cdots t_n^{\pm 1}]$ as
$B$ is invertible.  This way we can choose different co-ordinates which give
the same Laurent polynomial ring.

\paragraph*{Automorphism}  Any $B \in GL(n,\Z)$ defines an automorphism of
$\tilde{\tau}$ by
$$B X \otimes t^r =X \otimes t^{Br}$$
$$Bd (t^r)t^s =d (t^{Br}) t^{Bs}$$
Suppose $V$ is a module for $\tilde{\tau}$ then consider $\tilde{\tau}
\stackrel{B}{\to} \tilde{\tau} \to End V$.  By twisting by an automorphism we
get a new module.  This is equivalent to change of co-ordinates.
\paragraph*{Theorem} Suppose $\Omega_A/d_A$ acts as non-zero on $V$.  Then 
the following holds for a suitable choice of co-ordinates. 
\begin{enumerate}
\item[(1)] $k=n-1$ 
\item[(2)] If $t^r K_i \neq 0$ on $V$ then $i=n$ and  $r_n =0$ where $r=(r_1,
\cdots r_n) \in \Z^n$. 
\item[(3)] There exists non-zero central operators $z_1, \cdots z_{n-1}$ of
degree $(k_1, 0, \cdots 0)$, \\
$ \cdots , (0, \cdots 0, k_{n-1}, 0)$ where each $k_i$ is positive
integer. 
\item[(4)] There exists a proper submodule $W$ of $V$ for the Lie algebra 
$\tau \oplus \C d_n$ such that $V / W$ has finite dimensional weight space 
for the Cartan $h \oplus \sum \C K_i \oplus \C d_n$.  Further $t^r K_n$ acts as
scalar on $V/W$.
\end{enumerate}

\paragraph*{Corollary} Let $n \geq 2$.  Let $\overline{\tau}$ be the
quotient 
of $\tilde{\tau}$ by non-zero degree central operators.  Then
$\overline{\tau}$
does not admit representation with finite dimensional weight spaces where the center (degree zero) act non-trivially.
\paragraph*{Proof}  Here $k=0$.  Thus by the above theorem $\Omega_A/d_A$
acts trivially.  In particular each $K_i$ acts trivially.  

Let ${\cal G}_{aff} = {\cal G} 
\otimes \C [t_n, t_n^{-1}] \oplus \C K_n$ consider the Lie algebra 
homomorphism 

$$
\begin{array}{lll}
\varphi: \tilde{\tau} & \mapsto& {\cal G}_{aff} \otimes A_{n-1} \oplus D \\
X \otimes t^r & \mapsto & (X \otimes t_n^{r_n}) \otimes t^{r'},  
{\rm where} \  r'=(r_1, \cdots, r_{n-1}) \\
t^m K_i & \mapsto & \begin{cases} 0 \ {\rm if} \ 1 \leq i \leq n-1 \\
0 \ {\rm if} \  i=n, m_n \neq 0 \\
K_n \otimes t^m, \ {\rm if} \  m_n =0, i=n  
\end{cases} 
\\ [2mm]
d_i & \mapsto & d_i 
\end{array}
$$
From above theorem we see that on any irreducible integrable module for 
toroidal Lie algebra, the Ker $\varphi$ acts trivially upto the choice of 
co-ordinates.

Note that the quotient $V /W$ that occurs in the statement (4) of the theorem need not be 
irreducible nevertheless it admits an irreducible quotient.  It is proved
in
[E8] that from this irreducible quotient it is possible to recover the original module $V$.  Thus we have the following

\paragraph*{Remark}  The study of irreducible modules with finite
dimensional weight
spaces for toroidal Lie algebra $\tilde{\tau}$ is reduced to the study of
irreducible modules for ${\cal G}_{aff} \otimes A_{n-1} \oplus \C d_n$
with finite dimensional weight spaces where the infinite center $K_n
\otimes A_{n-1}$ acts as scalars.
\paragraph*{1.4}  Classification of irreducible integrable modules with
finite dimensional weight spaces [E4], [E5] and [Y].
\paragraph*{Definition} A weight module $V$ of $\tilde{\tau}$ is said
to be integrable if for every $v$ in $V$ and $\alpha \in \triangle, m \in \Z^n$ there exists $k = k (\alpha, m,v)$ such that
$$(X_{\alpha} \otimes t^m)^k v=0.$$
We will now describe the classification of irreducible integrable modules 
for the Toroidal Lie algebra.  

For each $i, 1 \leq i \leq n$, let $N_i$ be
a positive integer.  Let $\underline{a}_i = (a_{i1}, \cdots a_{iN_i})$ be 
non-zero distinct complex numbers.  Let $N=N_1 N_2 \cdots N_n$.  Let 
$I= (i_1, \cdots, i_n), 1 \leq i_j \leq N_j$ and  let $m= (m_1, \cdots,
m_n) 
\in \Z^n$. Define $a_I^m= a_{1i_1}^{m_1} \cdots a_{n i_n}^{m_n}$.

Let ${\cal G}_{aff} = {\cal G} \otimes \C [t_{n+1} t_{n+1}^{-1}] \otimes
\C K_{n+1}$.  Let $S= \{ I =(i_1, i_2, \cdots i_n ) \mid 1 \leq i_j
\leq N_j\}$ 
and  clearly $\#S=N$. Let $I_1, \cdots I_N$ be some order of
elements in $S$.  For each $i, 1 \leq i \leq N$ let $\lambda_i$ be a dominant 
integral weight for the affine Lie algebra ${\cal G}_{aff}$.  Let $V(\lambda_i)
$ be  the irreducible integrable highest weight module for ${\cal
G}_{aff}$. We will now define a ${\cal G}_{aff} \otimes A \oplus \C d_n$
module structure on
$$V (\lambda_1)\otimes \cdots \otimes V (\lambda_N) \otimes A \ {\rm by}$$
$$X \otimes t^r (v_1 \otimes \cdots \otimes v_N \otimes t^s)$$
$$= \displaystyle{\sum_{j=1}^{N}} a_{I_j}^{r} v_1 \otimes \cdots X v_j 
\otimes  \cdots \otimes v_N \otimes t^{s+r}$$
$$ d_i v_1 \otimes \cdots \otimes v_N \otimes t^s
=s_i(v_1 \otimes v_2 \otimes \cdots v_N \otimes t^s), 1 \leq i \leq n$$
$$d_{n+1} v_1 \otimes  \cdots \otimes v_N \otimes t^s$$
$$= \displaystyle{\sum_{j=1}^{N}} \ v_1 \otimes \cdots d_{n+1} v_j \otimes
\cdots \otimes v_N  \otimes t^s.$$ 
The following hold for the above module.  Let $h' = h_0 \oplus \sum \C K_i
\oplus 
\C d_i$.

\begin{enumerate}
\item[(1)] Integrable
\item[(2)] Finite dimensional weight spaces with respect to $h'$.
\item[(3)]  Completely reducible and most often irreducible.  All components
are isomorphic upto grade shift.
\end{enumerate}
\paragraph*{Theorem}  Any irreducible integrable module for the toroidal Lie 
algebra with finite dimensional weight spaces where center acts non-trivially
is isomorphic to a component of the above module up to a choice of co-ordinates.
\paragraph*{Remark}  The classification for the case where the center zero is
similar.  The ${\cal G}_{aff} $ has to be replaced by the finite 
dimension simple Lie algebra ${\cal G} $ and $V (\lambda_i)$ is 
irreducible integrable module for ${\cal G}$.  The case center acts trivially has also been done by Youngsun Yoon [Y].  She is able to give necessary and
sufficient condition for the above module to be irreducible.

\paragraph*{Remark}  It is well known that affine Kac-Moody Lie-algebras has
a rich representation theory like highest weight integrable modules where
Weyl character formula holds.  But the corresponding theory for
toroidal Lie algebras is very disappointing as the integrable modules are 
made out of affine modules.  It could still be possible  that semi integrable modules for toroidal Lie algebras could be of interest as in the case of affine Lie algebras [W].  It should be mentioned that integrable modules for
toroidal Lie algebras need not be completely reducible (even with the
assumption that each weight space is finite dimensional).  Such modules are
constructed by Chari and Le Thang [CL]  in two variables case.  So the
classification of integrable modules not necessarily irreducible is an 
interesting problem as their characters are $W$-invariant and may hold
modular properties.  There are also constructions of non-integrable modules
for toroidal Lie algebras by Ben Cox [C] and  Jing, Misra and Tan [JMT]
which may lead to Wakimoto type modules for toroidal Lie algebra.
\paragraph*{1.5}  Representations for toroidal Lie algebras through the
use of Vertex operators [EMY] and [EM].
  
Let $\Gamma $ be a non-degenerate integral lattice with symmetric bilinear form (,).  Let $h=\C \otimes_{\Z} \Gamma$.  Extend the form to $h$ bilinear.
For each $n$ let $h (n)$ be  a vector space isomorphic to $h$ through the 
isomorphism $a \to a (n), \ a \in h$.  We  create the Heisenberg Lie
algebra
$$A= \displaystyle{\oplus_{n \in \Z}}  h (n) \oplus \C C$$
with Lie structure
$[a (n), b(m)]= \delta_{m+n,0} \ n (a,b) C $ and $C$ is central.
Let ${\cal A}_{\pm} = \displaystyle{\oplus_{n \stackrel{>}{<}0}} h (n)$ so that 
$${\cal A} = {\cal A}_{+} \oplus h (0) \oplus \C C \oplus {\cal A}_{-}.$$
Let $C[\Gamma]$ be the group algebra and form the Fock space
$$V(\Gamma) = \C [\Gamma] \otimes S ({\cal A}_{-})$$
where $S({\cal A}_-)$ is the symmetric algebra of ${\cal A}_-, \ 
V(\Gamma)$
admits a natural Heisenberg Lie algebra structure.  Let $Q$ be even integral
sub lattice.

For $ \alpha \in Q$ let
$$X(\alpha, z) =z^{ \frac{(\alpha, \alpha)}{2}} \ exp \ T_- (\alpha, z)
e^{\alpha} z^{\alpha (0)} \ exp \ T_{+} (\alpha, z)$$
where $T_{\pm} (\alpha, z) = - \displaystyle{\sum_{n \stackrel{>}{<}0}} \
\frac{1}{n} \ \alpha (n) z^{-n}$ \\
$$z^{\alpha (0)} \cdot e^{\mu} \otimes u =z^{(\alpha, \mu)} e^{\mu} \otimes u$$
$$e^{\alpha} \cdot e^{\mu} \otimes u = e^{\alpha+ \mu} \otimes u.$$
There are all operators on 
$V[\Gamma] \otimes \C [z,z^{-1}]$.  Write $X (\alpha, z) =\displaystyle{\sum_{n \in \Z}} X_n (\alpha) z^{-n}$.  Then each component $X_n (\alpha)$ is an 
infinite sum but $X_n (\alpha) v$ is finite for any $v$ in $V[\Gamma]$.

We now choose simple finite dimensional Lie algebra ${\cal G}$ whose 
dynkin diagram is simply laced.  Now let $\Gamma$ be the lattice spanned by 
$\alpha_1,\cdots \alpha_d, \delta_1, \cdots \delta_n$, \\
$ \Lambda_1, \cdots \Lambda_n$ and $Q$ be the lattice spanned by $\alpha_1, 
\cdots \alpha_d, \delta_1, \cdots \delta_n$.
  We will now give a faithful and integrable representation for the toroidal Lie algebra in $n+1$ variables.

\paragraph*{Theorem} The following map defines a faithful representation for
the toroidal Lie algebra.
$$\tau \mapsto End V[\Gamma]$$
$$X_{\alpha} \otimes t^{(r, r_{n+1})} \mapsto X_{r_{n+1}} (\alpha+ \delta_r)$$
$$h \otimes t^{(r, r_{n+1})} \mapsto T_{r_{n+1}}^{ h} (\delta_r)$$
$$d(t^{(r, r_{n+1})})t^{(s,s_{n+1})} \mapsto \sum_{i=1}^{n} 
r_i T_{r_{n+1}+s_{n+1}}^{\delta_i} (\delta_{r+s})$$
$$+ r_{n+1} X_{r_{n+1}+s_{n+1} } (\delta_{r+s})$$
where $X_{\alpha}$ is a root vector of root $\alpha$ in ${\cal G}$.  $r=
(r_1, \cdots, r_n) \in \Z^n$ and
$$T_m^n (\gamma) = \displaystyle{\sum_{k \in \Z}} : \alpha (k) X_{n-k}
(\gamma): $$
\paragraph*{Remark ~ 1} The representation is highly reducible.  It is
proved in [E2] that there is a filtration such that successive quotients are loop modules.

\paragraph*{Remark ~ 2.}  The above construction of the representation has
been generalized to non-simply laced case by [JM1], [JY], [T2] and [XH].

Using this representation we can give a very nice presentation for the 
toroidal Lie algebra in two variables which has been used to construct modules
in [C] and [JMT].

\paragraph*{Theorem [EMY]} The following generators with the relations is a
presentation for the toroidal Lie algebra in two variables.

\paragraph*{ Generators}
$$C, \alpha_i^{\vee} (k), X_k (\pm \alpha_i), i=0, 1, 2 \cdots d, k \in
\Z$$

\paragraph*{Relations}
\begin{enumerate}
\item[(0)] $[C, \alpha_{i}^{\vee} (k)] = [C,X_k (\pm \alpha_i)] =0$
\item[(1)] $[\alpha_i^{\vee} (k), \alpha_j^{\vee}  (m)]= k  (\alpha_i, \alpha_j^{\vee}) \delta_{k+m,0} \ C$.
\item[(2)] $[\alpha_i^{\vee} (k), X_m (\pm \alpha_j)] = \pm (\alpha_i^{\vee}, 
\alpha_j) X_{m+k} (\pm \alpha_j)$.
\item[(3)] $[X_m (\alpha_i), X_n (-\alpha_j)]=- \delta_{ij} \{\alpha_i^{\vee}
(m+n) + 2m \delta_{m+n,0} \frac{C}{(\alpha_i, \alpha_j)}$.
\item[(4)] $[X_m (\alpha_i), X_n (\alpha_i)]= [X_m (- \alpha_i), 
X_n (-\alpha_i)]=0$.
$$
\left.
\begin{array}{ll}
ad (X_0 (\alpha_i))^{-A_{ji}} X_m (\alpha_j) &=0\\
(ad X_0 (- \alpha_i))^{A_{ji}} X_m (- \alpha_0) &= 0
\end{array}
\right \} \quad i \neq j
$$
\end{enumerate}

Here $A=(A_{ij})$ is an affine Cartan matrix of rank $d$.  $\{\alpha_0,
\cdots
\alpha_d\}$ is like set of simple roots and $\{ \alpha_0^{\vee}, \cdots \alpha_d^{\vee} \}$ is a set of simple co-roots and $(\alpha_i, \alpha_j)$ is the 
bilinear form given by the affine matrix $A$.  $(\alpha_i, \alpha_j^{\vee})$
is the evaluation of $\alpha_i$ at $\alpha_j^{\vee}$.
\paragraph*{$\S$ 2 ~~~  A generalization of Virasoro algebra}
 \paragraph*{2.1}  It is well known that the derivation algebra of $A= \C [t_1^{\pm 1}, \cdots t_n^{\pm 1}]$, denoted by Der$A$, is isomorphic to the
Lie algebra of diffeomorphisms of $n$ dimensional torus.  The universal
central extension of Der$A$ for $n=1$ is the well known Virasoro algebra.  The
representations of Virasoro algebra is studied in great detail.  See Ref.[KR].
It is interesting to note that Der$A, n \geq 2$ is centrally closed in the 
sense that it does not admit non-trivial central extension.  See 
Ref.[RSS]. But it admits an abelian extension (see [EM]) which generalizes
the Virasoro algebra.

In this section we recall the construction of Der$A$ modules due to Larsson 
[L3] and state its properties.  Later we construct modules for the abelian 
extension and state a conjecture on the classification of its modules.
\paragraph*{(2.2)}  Let $\C^n$ be  the $n$ copies of the complex field $\C$.  
Let $e_1, \cdots e_n$ be the standard basis of $\C^n$ and let $(,)$ be the 
standard form on $\C^n$ such that $(e_{i}, e_j)=\delta_{ij}$.  Let 
$\Z^n = \Z e_1 \oplus \Z e_2 \oplus \cdots \oplus \Z e_n$.  For $r \in \Z^n$ let $D^i (r)= t^r t_i \frac{d}{dt_i}$.  It is easy to verify that $D^i (r), 1 \leq i \leq n, r \in \Z^n$ is basis of Der$A$.  For $u \in \C^n$ let 
$$D(u,r) = \sum u_i D^i (r).$$ 
The Lie structure on Der$A$ is given by
$$[D(u,r), D(v,s)] =D (w,r+s)$$
where $w= (u,s) v - (v,r)u$.

Let $g \ell_n$ be the Lie algebra of $n \times n$ matrices with entries in 
$\C$.  Let $E_{ij}$ be the elementary matrix with $(i,j)$th entry to be 1
and zero elsewhere.  It is well known that $g \ell_n$ is spanned by 
$E_{ij}$ with Lie bracket
$$[E_{ij}, E_{k \ell}]= \delta_{jk} E_{i \ell} - \delta_{i \ell} E_{kj}.$$
Let $g \ell_n = s \ell_n \oplus \C I$ where $s \ell_n$ is the subalgebra of
trace zero matrices and $I$ is the identity matrix.  Let $V(\psi)$ be the 
irreducible finite dimensional module for $s \ell_n$ and let $I$ act by a scalar $b$.  Denote the resultant $g \ell_n$ module by $V(\psi, b)$.  For $\alpha \in \C^n$ we will make $F^{\alpha} (\psi,b)=V (\psi, b) \otimes A$ a Der$A$
module.  First denote $v \otimes t^m$ by $v(m) $ for $v$ in $V(\psi, b)$
and $m \in \Z^n$.  The following definition is due to Larsson [L3].
$$D(u,r)v(m)= (u,m+\alpha) v (m+r)+(\sum u_ir_j E_{ji}
v)(m+ r)$$
where $m,r \in \Z^n, u \in \C^n$ and $v \in V(\psi,b)$.
\paragraph*{Theorem [E3]}
\begin{enumerate}
\item[(1)]  $F^\alpha (\psi, b)$ is irreducible Der$A$ module if $(\psi,b)
\neq (\delta_{k},k), (0,b)$ where $\delta_k$ is the fundamental weight for
$s \ell_n, 1 \leq k \leq n-1$.
\item[(2)] $F^\alpha (0,b)$ is irreducible as Der$A$ module unless $\alpha
\in \Z^n$ and $b \in \{0,n\}$.
\item[(3)]  In all other cases $F^{\alpha} (\psi, b)$ is reducible and the
submodule structure has been worked out in [E3].

\end{enumerate}

The above result has been generalized by  Weiqiang Lin and Shaobin Tan in
[LT] for non-commuting variables, the so called quantum torus.

Let us note that $F^{\alpha}(\psi,b)$ is $A$ module by defining
$$t^r v (m) = v (m+r).$$

Then it is not difficult to see that $F^{\alpha} (\psi, b)$ is always
irreducible as \\
 $A \ltimes DerA$ module and the semi direct product is given by
$$[D(u,r), t^m ] = (u,m)t^{r+m}, \ [t^r,t^s]=0.$$
\paragraph*{Theorem [E6]}  Any irreducible $A \ltimes DerA$ module with
finite dimensional weight spaces has to be $F^{\alpha} (\psi,b)$ for some
$\psi$ dominant integral and some scalar $b$ and for some $\alpha$.
\paragraph*{(2.2)} We will now define an abelian extension for Der$A$
denoted by
$$\widehat{DerA} = \Omega_A/d_A \oplus DerA$$
\paragraph*{Definition} An abelian extension we mean an extension of
Der$A$ whose kernal is an abelian ideal.  

The Lie algebra structure on
$\widehat{DerA}$ is given in [EM].
Denote $K(u,r) = \sum u_i t^r K_i$
$$[D(u,r), D(v,s)] = D(w,r+s)- (u,s) (v,r) K(r,r+s)$$
where $w= (u,s)v- (v, r)u$.  
$$[D(u,r), K (v,s)]=(u,s)K (v,r+s)+(u,v) K(r,r+s)$$
$$[K (u,r), K(v,s)]=0.$$
We note that for $n=1, \  \widehat{DerA}$ is central extension and in fact
isomorphic to Virasoro algebra.
We will note that $GL(n,\Z)$ acts as automorphism on $\widehat{DerA}$.
Let
$B \in GL(n,\Z)$ which naturally acts on $\C^n$ and is denoted by dot.
$$B \cdot D (u,r) = D (B.u, B.r)$$
$$B \cdot K(u,r) = K ((B^T)^{-1} \cdot u,B \cdot r), \ u \in \C^n, r \in
\Z^n.$$
Applying $B$ on $\widehat{DerA}$ is what we call change of co-ordinates.

We will now construct modules for $\widehat{DerA}$ and we will work with
$(n+1)$ variables and $(n+1)$th variable has been chosen as $t_0$.
We will first define a $\Z$ grading on $\widehat{DerA}$ by declaring the
grade of $D(u,r)$ and $K(u,r)$ by $r_0$ where $r= (r_0, r_1, \cdots r_n)
\in \Z^{n+1}$.  Now $\widehat{DerA}_+ $ (respectively $\widehat{DerA}_-,
\widehat{Der A}_0)$ is the linear span of $D(u,r) $ and $K(u,r)$ for 
$u \in \C^{n+1}$ and $r_0 >0$ (respectively $r_0<0$ and $r_0 =0)$.  Then
clearly we have the following decomposition of subalgebras
$$\widehat{DerA} = \widehat{DerA}_- \oplus \widehat{DerA}_0 \oplus
\widehat{DerA}_+.$$
We will now make $V(\psi,b) \otimes A_n$ as $\widehat{DerA}_0$ module by
declaring
$$
\begin{array}{lll}
t^r t_0 \frac{d}{d0} \cdot v \otimes t^s &=&d  (v \otimes t^{r+s}) \\
t^r K_p \cdot v \otimes t^s &=&0, 1 \leq p \leq n \\
t^r K_0 \cdot v \otimes t^s &=& C (v \otimes t^{r+s}) 
\end{array}
$$
for scalars $d$ and $C$ and $r,s \in \Z^n$.  Let $\widehat{DerA}_+$ act
trivially  on $V(\psi,b) \otimes A_n$.  Now consider the induced module
for $DerA$.
$$J (\psi, b) =U (\widehat{DerA}) \otimes_{( \widehat{DerA}_0  \oplus 
\widehat{DerA}_+)} (V ( \psi, b) 
\otimes A_n).$$ 

Then $J (\psi, b)$ has a unique irreducible quotient (we assume either $C
\neq 0$ on $d \neq 0)$ say $I (\psi, b)$.  Now from Berman and Billig [BB]
it follows that $I (\psi,b)$ has finite dimensional graded spaces as $\Z^{n+1}$
graded vector space.

\paragraph*{Conjecture}  Any irreducible module for $\widehat{DerA}$ with 
finite dimensional weight spaces has to come from Larsson's construction or
a highest weight module in the above  sense upto suitable choice of 
co-ordinates.
\paragraph*{$\S$ 3  Full Toroidal Lie algebra}
 
In this section we define a more general toroidal Lie  algebra which we call
full toroidal Lie algebra and denoted by
$$\widehat{\tau} = {\cal G} \otimes A \oplus \Omega_A/d_A \oplus DerA.$$
The Lie structure for the first two components and for the last two 
component has already been defined.  Thus we defined Lie bracket between the first and the third component.
$$[D(u,r), X \otimes t^s] = (u,s) X \otimes t^{r+s}.$$
We will now state the results on classification of integrable modules for 
$\widehat{\tau}$.  As earlier a module for $\widehat{\tau}$ is integrable
if it is integrable as $\tau$ module.  That is $X_{\alpha} \otimes t^r$ acts as locally nilpotent on the module.  We will work with $(n+1)$ variables and
the $(n+1)$th variable is chosen as $t_0$.

Let us note that the span of $\{K_0, K_1, \cdots, K_n\}$ is the center of $\widehat{\tau}$ and on an irreducible integrable weight module each $K_i$
 acts as integer.  We will first construct irreducible integrable modules for
$\widehat{\tau}$ with finite dimensional weight  spaces where $K_0$ acts
as positive integer and $K_i(1 \leq i \leq n)$ acts trivially.  Clearly
$\widehat{\tau}$ is $\Z^{n+1}$ graded and we will extract one $\Z$-
gradation that is by $t_0$.  Let $\widehat{\tau}_+ $ (respectively 
$\widehat{\tau}_-, \widehat{\tau}_0)$ be the linear span of ${\cal G} 
\otimes t_0^{r_0}t^r, D(u, (r_0,r)), K(u, (r_0,r)) $ where $r_0 >0$
(respectively $r_0<0$
and  $r_0=0\}$.

Then $\widehat{\tau} = \widehat{\tau}_- \oplus \widehat{\tau}_0 \oplus \widehat{\tau}_+$.

We will first define $\widehat{\tau}_0$ module and then make $\widehat{\tau}_+$
acts trivially  and then we consider the standard induced highest weight module.  Let $A_n =\C [t_1^{\pm 1}, \cdots t_n^{\pm 1} ]$ and recall that
$F^{\alpha} (\psi, b) = V (\psi,b) \otimes A_n$ is a Der$A_n$-module.  Let
$W$ be finite dimensional irreducible module for the finite dimensional
 Lie algebra ${\cal G}$.  We fix a positive integer $C_0$ and a complex number $d$.

Consider $T= W \otimes F^{\alpha} (\psi, b)$ which will now  be made
into $\widehat{\tau}_0$-module.

Let $w \in W, v \in V (\psi, b)$ and $r,s \in \Z^n$.  \\
$X \otimes t^r \cdot w \otimes v \otimes t^s = X w \otimes v 
\otimes t^{r +s}, X \in {\cal G}$ \\
$D(u,r) w \otimes v \otimes t^s= w \otimes D(u,r) (v \otimes t^s), u \in
\C^n$,
$$D(e_0, r) w \otimes v \otimes t^s =d (w \otimes v \otimes t^{r+s})$$
$$K(u,r) w \otimes v \otimes t^s= 0 \ {\rm for} \ u \in \C^n$$
$$K(e_0,r) w \otimes v \otimes t^s =C_0 (w \otimes v \otimes t^{s+r}).$$

Recall $K(u,r)$ and $D(u,r)$ are linear in $u$.  We first defined $K(u,r)$ and
$D(u,r)$ for the last $n$ variable and then separately for $K(e_0,r)$ and
$D(e_0,r)$.  It is easy to check that $T$ is a module for
$\widehat{\tau}_0$ and in fact it is irreducible.  Let $\widehat{\tau}_+$ act trivially and consider the induced $\widehat{\tau}$-module 
$$M(T) = U (\widehat{\tau}) \otimes_{\widehat{\tau}_+ \oplus
\widehat{\tau}_0} T.$$

It is easy to see that $M(T)$ has a unique irreducible quotient say 
$V(T)$.  $V(T)$ has a natural $\Z^{n+1}$-gradation and it follows from Berman and Billig [BB] that each weight space is finite dimensional.  From the fact
that $W$ is finite dimensional it follows that $V(T)$ is integrable but that
is little non trivial.

\paragraph*{Theorem} ([JM2], [EJ]).  Suppose $V'$ is irreducible
integrable module for $\widehat{\tau}$ with finite dimensional weight spaces.  Suppose the 
center acts non-trivially.  Then $V'$  is isomorphic to $V(T)$ for some 
$T$ up to a suitable choice of co-ordinates.

The proof consists of three steps.  The first step is to prove that $V'$ is 
a highest weight module upto suitable choice of co-ordinates, which follows 
from [E5].  The second step is to reduce the classification of irreducible 
integrable highest weight module for the classification of $A_n \ltimes 
DerA_n$  irreducible modules whose proof can be found in [JM2].  The third
step is to classify 
$A \ltimes DerA$ modules which is done in [E6].

We will now describe the irreducible integrable modules for $\widehat{\tau}$
where the center acts trivially.

We will work with $n$ variables.  Let $W$ be any non-trivial finite 
dimensional irreducible module for ${\cal G}$.  Then we make 
$W \otimes V(\psi, b) \otimes A$  as $\widehat{\tau}$-module.
$$X \otimes t^r \cdot u \otimes v \otimes t^s =X u \otimes v \otimes t^{r+s}$$
$$D(u,r) u \otimes  v \otimes t^s = u \otimes D(u,r) (v \otimes t^s) $$
$$K(u,r) u \otimes v \otimes  t^s =0.$$

It is easy to see that $W \otimes V (\psi,b) \otimes A$
is an irreducible integrable module for $\widehat{\tau}
$ with finite dimensional weight spaces.  Note that $W$ has to be 
non-trivial for the module to be irreducible.

\paragraph*{Theorem}  ([JM2], [EJ]) ~~ Suppose $V'$ is an irreducible integrable module for $\widehat{\tau}$ with finite dimensional weight space.  Suppose
 the center acts trivially and suppose ${\cal G}$   action on $V'$ is 
non-trivial.  Then $V' \cong W \otimes V(\psi,b) \otimes A $ for some
$W$ non-trivial module for ${\cal G}$  and for some $\psi$ and $b$.

\paragraph*{Remark}  The assumption that ${\cal G}$ acts non-trivially
is necessary in the theorem.

The case where ${\cal G}$ acts trivially is open.  In this case the
problem is nothing but the conjecture stated in section 2.

It is an open problem for a long time to explicitly construct a module 
for $\widehat{\tau}$.  Several attempts have been made in ([EM], [E7], [BB],
[BBS]) to find a realization for $\widehat{\tau}$.  Eventually in a  
remarkable paper [B3], Yuly Billig succeeded in construction of  a module
for $\tilde{\tau}$ through the use of Vertex Operator Algebras.

\paragraph*{Remark}  Consider $D^* =\{D(u,r) \in DerA \mid (u,r)=0\}$.
Then
${\cal G} \otimes A \oplus  \Omega_A /d_A \oplus D^*$ is a subalgebra of 
$\widehat{\tau}$ and form an important class of examples of the so called 
Extended Affine Lie Algebra (EALA).

\paragraph*{Definition of EALA}

Let  $L$ be a Lie algebra over $\C$.  Assume that \\
(EA1)  $L$  has a nondegenerate invariant symmetric bilinear form (,) \\
(EA2) $L$  has a nontrivial finite dimensional self-centralising 
ad-diagonalisable abelian subalgebra ${ H}$.

We will assume 3 further axioms about the triple $(L, (\cdot , \cdot), H)$.
To describe these axioms we need further notation.  Using (EA2), we have

$$L = \displaystyle{\bigoplus_{\alpha \in H}} \ L_{\alpha},$$
and
$$L_0 =H.$$
where $L_{\alpha} =\{x \in L: [h,x] = \alpha (h) x \ {\rm for \ all} \ h
\in H\}$, \ and \ $H^* $ is the complex dual space of $H$.  Let
$$R= \{\alpha \in H^* : L_{\alpha} \neq \{0\}\}.$$
$R$ is called the {\it root system of} $L$.  Note that since $H \neq \{0\}$
we have $0 \in R$.  Also,
$$\alpha, \beta \in R, \alpha+ \beta \neq 0 \Rightarrow (L_{\alpha}, L_{\beta})
= \{0\} \eqno{1.1}$$
Thus, $-R=R$.  Moreover, $(\cdot , \cdot)$ is nondegenerate on $H$.  So as usual we can transfer $(\cdot,\cdot)$ to a form on $H^*$.  Let 
$$R^* =\{ \alpha \in R: (\alpha, \alpha) \neq 0\} \ {\rm and} \ R^0 =
\{\alpha \in R: (\alpha, \alpha) = 0\}.$$
The elements of $R^*$ (respectively $R^0)$ are called non isotropic (respectively isotropic ) roots.  We have
$$R=R^0 \cup R^*.$$
(EA3) $\alpha \in R^*, x_{\alpha} \in L_{\alpha} \Rightarrow adx_{\alpha}$
acts locally nil-potently on $L$. \\
(EA4) $R$ is a discrete subset of $H^*$. \\
(EA5) $R$ is irreducible.  That is

(a) $R^* =R_1 \cup R_2, (R_1,R_2) =\{0\} \Rightarrow R_1= 0 \ {\rm or} \
R_2=0$ \\

(b) $\sigma \in R^0 \Rightarrow$. There exists $\alpha \in R^*$ such
that $\alpha+ \sigma \in R$.

Then $L$ is called EALA.  

Extensive research has been done on the structure and
classification of EALA's.  See [AABGP] and [AG] and the references there in.

\paragraph*{Open Problems.}  Certainty toroidal Lie algebras form an 
important class of examples of Extended Affine Lie Algebras (EALA)  and 
certain progress has been made towards representation theory of toroidal 
Lie algebras.  One question is that can we develop similar theory for EALA.  
In particular  is it possible to classify irreducible integrable modules 
for EALA's.

Apart from toroidal Lie algebras, one particular example of EALA for which 
representations are studied is an EALA whose co-ordinated algebra is a 
quantum torus and they generally appear in type A.  The representations of
this Lie algebra is studied by Yun Gao, Kei Miki and myself.  It will
be good idea to classify irreducible integrable modules for this Lie algebra.

The classification of integrable modules for the toroidal super-algebra has
been attempted in [EZ].  Several open problems exists in super case.

More complex algebras, the so called quantum toroidal is also pursued by several authors like E. Vasserot, M. Varagnolo, Kei Miki, Naihuan Jing and Oliver Schiffmann.

\paragraph*{Applications.}  The vertex representation for toroidal Lie 
algebras for the homogeneous picture  is done in [EM].  Similar construction
is made for the principal picture in [B1].  Both constructions have found
applications in differential equations via [B2], [ISW] and [IT].

The papers [EM] and [EMY] have been mentioned in physics literature from
the works of T.A. Larsson, M. Calixto,, T. Iami, T.Ueno and H. Kanno.  
In fact Larsson has reinterpreted the results of [EM] and [EMY] in the
language of Physics in [L4] and [L5].

\pagebreak

\begin{center}
{\bf REFERENCES}
\end{center}
\vskip 7mm
\begin{enumerate}
\item[{[AABGP]}]  Allison, B.N., Azam, A., Berman, S., Gao, Y. and
Pianzola, A., Extended Affine Lie algebras and their root systems,
Mem. Amer. Math. Soc. 126 (1997), No.603, 1-122.
\item[{[AG]}] Allison, B.N. and Gao Y. The root system and the core  of an 
extended affine Lie algebra,  Selecta Math. (N.S.), 7 (2001), 
149-212.
\item[{[Ba]}] Batra, P. Representations of twisted multi-loop algebras, 
J. Algebra, 272 (2004), 404-416.
\item[{[BB]}]  Berman, S. and Billig, Y. Irreducible representations for toroidal Lie algebras, J. Algebra, 221 (1999), 188-231.
\item[{[BBS]}] Berman, S., Billig, Y. and Szmigielski, J. Vertex Operator
algebras and the representation theory of toroidal Lie algebras, Contemp.
Math. 297 (2002), 1-26.
\item[{[BC]}] Berman, S. and Cox, B. Enveloping algebras and representations
on toroidal Lie algebras, Pacific J. Math. 165 (1994), 239-267.
\item[{[B1]}]  Billig, Y. Principal Vertex operator representations for toroidal Lie algebras, J. Math. Phys. 39 (1998), No.7, 3844-3684.
\item[{[B2]}] Billig, Y. An extension of the Korteweg-de  Vries hierarchy arising from a representation of toroidal Lie algebra, J. Algebra, 217
(1999), 40-64.
\item[{[B3]}]  Billig, Y. Energy-momentum tensor for the toroidal Lie
algebra, Arxiv. Math.RT/0201313.
\item[{[CL]}] Chari, V. and Le Thang, Representations of double affine Lie
algebras, A tribute to C.S. Seshadri, Trends Math., 199-219, Birkhauser,
(2003).
\item[{[C]}] Cox Ben. Two realizations of toroidal $s \ell_2 (\C)$,
Contemp. Math. 297 (2002), 47-68.
\item[{[CF]}] Cox B and Futorny, V. Borel Subalgebras and categories of highest weight modules for toroidal Lie algebras, J. Algebra, 236 (2001), 
1-28.
\item[{[DFP]}] Dimitrov, I., Futorny, V. and Penkov, I. A reduction theorem
for highest weight modules over toroidal Lie algebra, Comm. Math. Phys.
250 (2004), 47-63.
\item[{[DMS]}] Di Francesco, P., Mathieu, P. and Senechal, D. Conformal 
Field Theory (springer, New York, 1997).
\item[{[E1]}] Eswara Rao, S. Representations  of Witt algebras, Publ.
Res. Inst. Math. Sci., 30 (1994), 191-201.
\item[{[E2]}] Eswara Rao, S. Iterated loop modules and a filtration for
vertex  representation of toroidal Lie algebra, Pacific J. Math.
171(1995), 511-528.
\item[{[E3]}] Eswara Rao, S. Irreducible representations of the Lie algebra of the diffeomorphisms of a $d$-dimensional torus, J. Algebra 
182 (1996),401-421.

\item[{[E4]}] Eswara Rao, S. Classification of irreducible integrable modules
for multi-loop algebras with finite dimensional weight spaces, J. Algebra 246 (2001), 215-225.
\item[{[E5]}]  Eswara Rao, S. Classification of irreducible integrable modules for toroidal Lie algebras with finite dimensional weight spaces, J. 
Algebra, 277 (2004), 318-348.
\item[{[E6]}] Eswara Rao, S. Partial classification of modules for Lie 
algebra of diffeomorphisms of $d$-dimensional torus, J. Math. Phys. 
45 (2004), No.8, 3322-3333.

\item[{[E7]}] Eswara Rao, S. Generalized Virasoro operators, Comm. Algebra, 
32 (2004), No.9, 3581 - 3608.
\item[{[E8]}] Eswara Rao, S. Irreducible representations for toroidal Lie 
algebras, Arxiv. Math. RT/0212137, To appear in J. Pure Appl. Algebra.
\item[{[EJ]}] Eswara Rao, S. and Jiang, C. Classification of irreducible 
integrable representations for the full toroidal Lie algebras,  Preprint (2004).
\item[{[EM]}] Eswara Rao, S and Moody, R.V. Vertex representations for
 $n$-toroidal Lie algebras and a generalization of the Virasoro algebra, Comm.
 Math. Phys. 159 (1994), 239-264.
\item[{[EMY]}]  Eswara Rao, S. Moody, R.V. and Yokonuma, T. Toroidal Lie
algebras and Vertex representations, Geom. Dedicata, 35 (1990), 283-307.
\item[{[EZ]}] Eswara Rao, S and Zhao, K. On integrable representations of 
toroidal Lie super-algebras, Contemp. Math. 343 (2004), 243-261.
\item[{[FM]}] Fabbri, M. and Moody, R.V. Irreducible representations of
Virasoro - toroidal Lie algebras, Comm. Math. Phys, 159 (1994), 1-13.
\item[{[FR]}]  Figueirido, F and Ramos, E. Fock space representation of
the algebra of diffeomorphisms of the $n$-torus, Phys. Lett. B. 238 (1990), 247-251.
\item[{[FK]}]  Frenkel, I. and Kac, V.G., Basic Representations of
affine Lie algebras and dual resonance models, 
Invent. Math. 62 (1980), 23-66.
\item[{[FK]}] Futorny, V and Kashuba, I. Verma type modules for toroidal Lie algebras, Comm. Algebra 27 (1999), No.8, 3979-3991.
\item[{[IT]}] Ikeda, T. and Takasaki, K. Toroidal Lie algebras and Bogoyavlensky's $2+1$ dimensional equation, Int. Math. Res. Not. 2001, No.7,
329-369.
\item[{[ISW]}]  Iohara, K. Saito, Y and Wakimoto, M. Hiroto bilinear forms with
2-toroidal symmetry, Phys. Lett. A. 254 (1999) No.1-2, 37-46.
\item[{[JM1]}] Jiang, C. and Meng, D. Vertex representations of $n+1$ toroidal Lie algebra of type B.  J. Algebra 246 (2001), 564-593.
\item[{[JM2]}] Jiang, C and Meng, D. Integrable representations for 
generalized Virasoro-toroidal Lie algebra, J. Algebra, 270 (2003),
307-334.
\item[{[JY]}] Jiang, C and Hong You. Vertex operator construction for
the toroidal lie algebra of type $F_4$, Comm. Algebra, 31 (2003), No.5,
2161-2182.
\item[{[JMT]}] Jing, N. Misra K.C. and Tan, S. Bosonic realizations of higher level toroidal Lie algebras, to appear in Pacific J. Math. (2004).
\item[{[K]}] Kac, V. Infinite dimensional Lie algebras, 3rd edition, Cambridge
University Press, Cambridge, U.K. 1990.
\item[{[Ka]}] Kassel, C., Kahler differentials and coverings of complex simple Lie Algebra extended over commutative algebras, J. Pure Appl. Algebra 34 (1985), 266-275.
\item[{[KR]}]  Kac, V and Raina, A.K. Bombay lectures on highest weight
representations of infinite dimensional Lie algebras, World Scientific, Singapore 1987.
\item[{[L1]}]  Larsson, T.A. Multidimensional Virasoro algebra, Phys.
Lett. B 231 (1989), 94-96.
\item[{[L2]}]  Larsson, T.A. Central and non-central extensions of multi-
grade Lie algebras, J. Phys. A 25 (1992), 1177-1184.
\item[{[L3]}] Larsson, T.A. Conformal fields: A class of representations
of
vect (N). Internat. J. Modern Phys. A  7 (1992), 6493-6508.
\item[{[L4]}] Larsson, T.A. Lowest energy representations of non-centrally
extended diffeomorphism algebras, Comm. Math. Phys. 201 (1999), 461-470.
\item[{[L5]}] Larsson, T.A., Extended diffeomorphism algebra and
trajectories in Jet space, Comm. Math. Phys, 214 (2000), 469-491.
\item[{[LL]}] Lepowsky, J. and Li, H., Introduction to vertex operator algebras and their representations, Progress in Mathematics, 227, Birkhäuser, Boston, 2004.
\item[{[LT]}] Lin Weiqiang and Tan Shaobin, Representations of the Lie algebra of derivations for quantum torus, J. Algebra, 275 (2004), 250-274.
\item[{[MY]}]  Morita, J. and Yoshi, Y. Universal central extensions of
Chevally algebras over Laurent polynomial rings and GIM Lie algebras, Proc. Japan Acad. Ser.A, 61 (1985), 179-181.
\item[{[RS]}] Ramos, E and Shrock, R.E. Infinite dimensional $\Z^n$-indexed Lie algebra and their super symmetric generalizations, Internat. J. Modern Phys. A 4 (1989), 4295-4302.
\item[{[ RSS]}] Ramos, E., Sah, C.H. and Shrock, R.E. Algebra of diffeomorphisms of the $N$-torus, J. Math. Phys., 31 (1990) 1805-1816.
\item[{[T1]}]  Tan Shaobin, Vertex operator representations for toroidal Lie algebra of type $A_1$, Math. Z. 230 (1999)  No.4, 621-657.
\item[{[T2]}]  Tan Shaobin, Vertex operator representations for toroidal Lie algebra of type $B$, Comm. Algebra 27 (1999), No.8, 3593-3618.
\item[{[XH]}] Xia, L.M. and Hu, N.H., Irreducible representations for Virasoro-
toroidal Lie algebra,  J. Pure Appl. Algebra 194 (2004), 
211-237. 
\item[{[Y]}] Youngsun Yoon, On polynomial representations of current 
algebras, J. Algebra 252(2) (2002), 376-393.
\item[{[W]}] Wakimoto, M., Lectures on infinite dimensional Lie algebras, 
World Scientific Publishing Company, NJ (2001).
\end{enumerate}

\end{document}